\documentclass[12pt]{article}

\usepackage{amssymb}
\usepackage{amsmath}
\usepackage{amsfonts}
\usepackage{epsfig}
\usepackage{color}
\usepackage{float}

\newtheorem{theorem}{Theorem}[section]
\newtheorem{prop}[theorem]{Proposition}

\newtheorem{lemma}[theorem]{Lemma}

\newtheorem{rem}[theorem]{Remark}

\newtheorem{rems}[theorem]{Remarks}

\newcommand{\sn}{\smallskip\noindent}
\newcommand{\ha}{\textstyle{\frac{1}{2}}\displaystyle}

\title{\bf Envelopes of lines, unfoldings and breaking symmetry}

\date{}
\author{Peter Giblin \and Alexander Wettig}

\begin{document}

\maketitle

\section{Introduction}

\begin{figure}[!h]
\begin{center}
\includegraphics[width=2in]{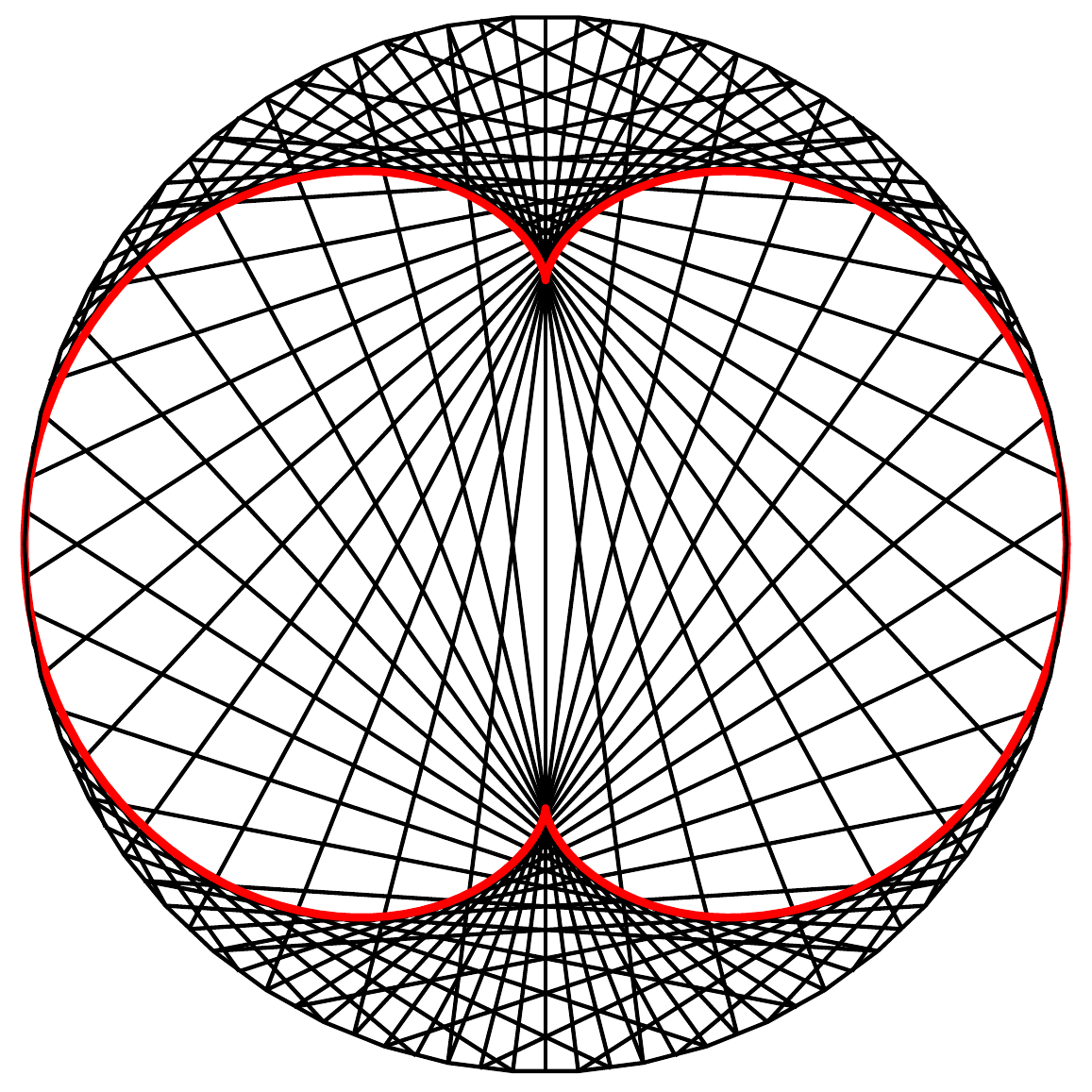}
\includegraphics[width=2in]{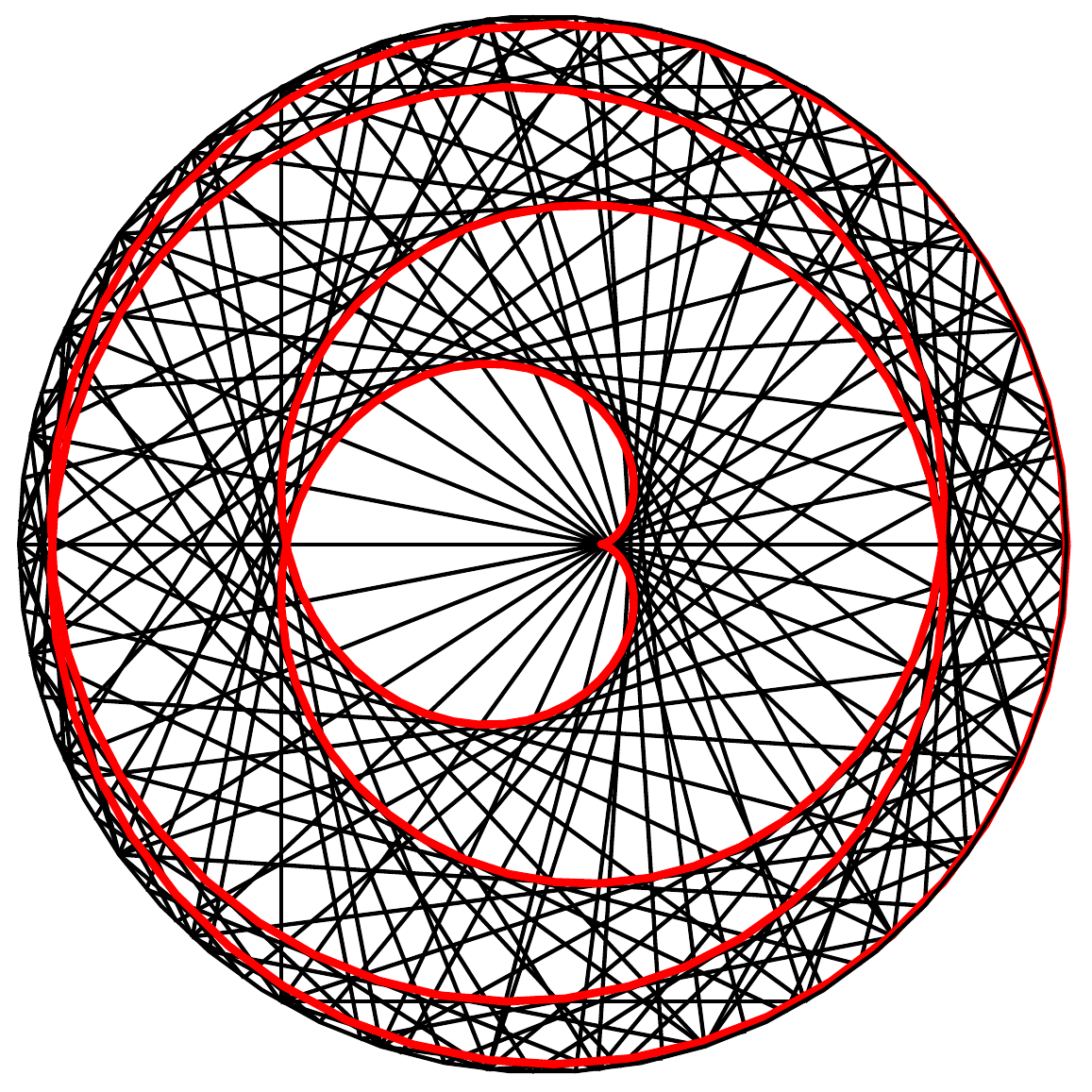}
\includegraphics[width=2in]{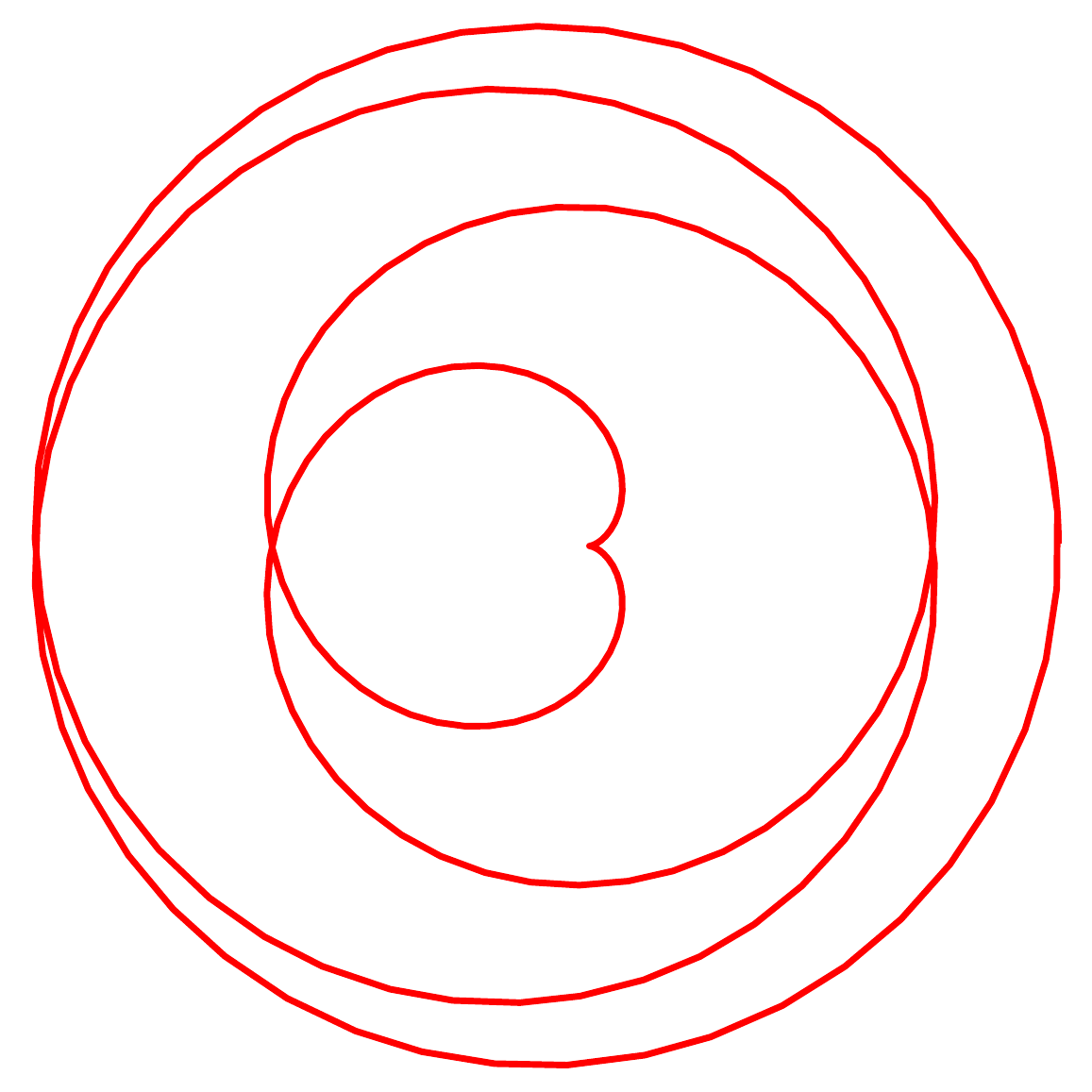}
\end{center}
\caption{\small All the figures show the ``embroidery'' envelope and also the circle.
Left:  $m=3$, with 2 cusps and no self-crossings.  Centre: $m=\frac{5}{4}$, which has 1 cusp and 3 self-crossings.
Right: just the envelope of the center diagram without the lines being drawn.}
\label{fig:exs1}
\end{figure}

It is very striking that when certain families of straight lines are drawn in the plane
it is not only the lines which are apparent to the eye but also a {\em curve} which has
the lines for its tangents. A simple example is shown in Figure~\ref{fig:exs1}, left where
the curve, called the {\em envelope} of the lines, has two lobes which join at
two ``sharp points'' called {\em cusps}. The center diagram of this figure is
more complicated: the curve drawn heavily is in fact tangent to all the lines and it
has one cusp (and also the curve crosses itself in three places, all on the
horizontal axis of symmetry). To make this more obvious the envelope alone is
drawn in the right-hand diagram. These diagrams of lines and their envelope
are often called ``embroidery diagrams'' since the lines are all formed by joining two points
of a unit circle, in rather the way that cotton might be  stretched across a circular
embroidery hoop.  For each envelope 
we choose a rational number $m=a/b$, where $a$ and
$b$ are coprime, then join the point $(\cos t, \sin t)$ to the point $(\cos mt, \sin mt)$ for
many values of $t$. The respective values of $m$ are given in the caption. We show how
to count the cusps and tangencies with the circle in \S\ref{s:1circle}. 
The value of $m$ will be restricted to {\em rational} values to ensure a closed
envelope curve: there is no continuously varying parameter in this construction which
allows us to observe continuous changes or evolutions of the envelope curve.

In this article we briefly explore these envelopes in \S\ref{s:1circle}, but then move on to some
generalizations, also involving envelopes of straight lines, but with continuous
parameters which can be varied. Envelopes constructed from straight lines
 are a natural and simple
way in which to construct curves with {\em singularities}
or {\em singular points}, that is points at which
the speed of the curve becomes zero;  usually this manifests itself a a cusp (sharp
point) such as those in Figure~\ref{fig:exs1}. There are other circumstances in which
``higher cusps'' occur and it is the way in which these break up and evolve---the
technical and rather descriptive word is ``unfold''---as the parameters change which will be our main topic. 
Such unfoldings are a major concern of ``Singularity Theory'' and this article is
initended to introduce some of the ideas in a concrete context.

In \S\ref{s:2circles} we introduce one continuous parameter $r>1$:
the straight lines join corresponding points $(\cos t, \sin t)$ and
$(r\cos mt, r\sin mt)$ of two circles centered at the
origin, of radii
1 and $r$.  (Here $m$ is a rational number which, in any one family of
envelopes, will be held fixed.) In
\S\ref{s:higher} we allow the center of the second circle to move  from $(0,0)$
to $(0,d)$, introducing a second continuous parameter $d$: we now have a two-parameter
family of lines with parameters $r,d$. The point of
doing this is to show how singularities evolve as the
two parameters $r, d$ change. As above, this process is 
called ``unfolding the singularity'', 
 revealing its inner structure. Some unfoldings may not be sufficient to
 reveal the full structure; those that do are called ``versal unfoldings''.
 There is an illustration in Figure~\ref{fig:butterfly1}, which shows how the ``butterfly singularity'', which is the ``sharp point'' to the left in the center diagram, unfolds completely (versally) by
varying the parameters $r,d$ with $m$ fixed at $3/4$.  Note that until we see the unfolding
there is no particular visual difference between the two sharp points in the center
diagram, but the one further to the right appears unaffected by the unfolding---as indeed
it is. This cusp is stable, there is nothing further to reveal by unfolding.

 Unfolding is one of the central concepts of
Singularity Theory and Catastrophe Theory, 
originally introduced by Ren\'{e} Thom and Hassler Whitney
in the 1950s and 1960s and then hugely developed by many others
including John Mather. There is some information on this in \cite{W,Mon,Wiki}
but no comprehensive history of the subject as yet exists.

Using readily available online graphics programs it is possible to explore these
envelopes experimentally and to witness the changes in the number and nature of
the singularities as parameters vary. Writing the instructions for such programs
requires calculation of the envelope curves starting from the family of lines and we
include some standard results and calculations
 in this direction in \S\ref{s:general}. The two-parameterl case
studied in \S\ref{s:higher} is available in the Desmos interactive program \cite{Desmos};
in \S\ref{s:general} and elsewhere
we shall explain where the formulas in this program come from. The program allows the
rational number $m=a/b$ to be chosen, the
lines and the envelope curve to be drawn or withheld and the parameters
$r =$ the radius of the second circle and $d=$ the displacement of the center 
in the $y$-direction to be changed. In fact the center of the second circle
can be moved more generally to $(c,d)$. 

In Remark \ref{rem:caustic} we explain briefly how some of the envelopes in \S\ref{s:2circles}  can be generated in a different way.

\medskip

The envelopes considered in this article are examples of a wide class where 
corresponding points of two
plane curves $C_1$ and $C_2$ are joined by straight lines. The two curves can
also coincide, and the ``correspondence'' can be any
relationship, differentiable when expressed in local parameters; it does not have to be one-to-one. For example the envelope of lines joining
points of  $C_1$ and $C_2$ at which the tangents are parallel gives the {\em center symmetry
set}, initially studied by Janeczko \cite{SJ} but later developed into a substantial
theory extending to higher dimensions by many authors; see \cite{GZ} for example.
In \cite{DZ} the authors study the related ``Wigner caustic'' which is the envelope of lines {\em halfway between} parallel tangents and which first arose in physics. In all cases a central
question is the nature of the singular points (cusps etc.) and the extent to which they
are fully unfolded, with their inner structure revealed, by whatever 
smoothly varying parameters are
available.

\section{An envelope from one circle}\label{s:1circle}

A well-known way of forming an envelope from chords of a single circle is  as follows (see for example \cite[\S5.7(2), \S7.14(6)]{CS}). For each point $(\cos t, \sin t), \ 0\le t \le 2\pi$, on the unit
circle construct the chord from this point to $(\cos(mt), \sin(mt))$ where $m$ is a positive integer $>1$. (For $t=0$ we use the tangent at the origin instead
of a chord.) These lines have an envelope which
has $m-1$ cusps and no self-crossings. See Figure~\ref{fig:exs1}, left, for the case $m=3$.  To generalise, we let
$m=\frac{a}{b}$ be any rational number, not equal to 0 or $\pm 1$ and in its lowest terms. We consider three problems, in increasing order of difficulty: (i) find the range of values of
$t$ which allows the envelope to close; (ii) find the number of cusps on the envelope; and (iii) find the number of self-crossings on the envelope.
Two examples are given also in Figure~\ref{fig:exs1}.

We shall prove the following.
\begin{theorem}
For the envelope constructed as above, with $m \ne -1,0,1$ a rational number $\frac{a}{b}$ in its lowest terms, with $b>0$,

{\rm (i)} the envelope, starting at $t=0$,  closes for $t=2b\pi$,

{\rm (ii)} the number of cusps is $|a-b|$, occurring when $(m-1)t$ is an odd multiple of $\pi$; there are also $|a-b|$ points of tangency between the envelope and the circle,
occurring when $(m-1)t$ is an {\em even} multiple of $\pi$.

In addition it can be shown that

{\rm (iii)} the number of self-crossings is
\begin{eqnarray*}
(|a|-1)|a-b|  &\mbox{if} & |m| <1 \\
(b-1)|a-b|  &\mbox{if} & |m| >1.
\end{eqnarray*}
\label{th}
\end{theorem}
As examples, see those in Figure~\ref{fig:exs1}.
See \S\ref{s:general} below for a more general discussion of envelopes of lines.

First, we shall need a formula for the envelope itself.  The line joining the points $(\cos t, \sin t)$ and $(\cos mt, \sin mt)$ ($m\ne 1)$ has equation
\begin{equation}
 x(\sin mt-\sin t) - y(\cos mt -\cos t) = \sin(m-1)t.
\label{eq:F1circle}
\end{equation}
The standard way to obtain the envelope of such a family of lines---that is, a curve which is tangent to all the lines---is to differentiate with respect to $t$ and then to solve for $x$ and $y$
as functions of $t$ from the two equations.  The solution will be written
as $x=X(t), y=Y(t)$ in what follows. After some trigonometric work we find
\begin{equation}
X(t) = \textstyle{\frac{1}{m+1}}\displaystyle(m\cos t + \cos mt), \ \ Y(t)=\textstyle{\frac{1}{m+1}}\displaystyle(m\sin t + \sin mt).
\label{eq:xy}
\end{equation}

\subsection{Tangencies with the circle}
It is important to note that in the process of deriving this parametrization there is a cancellation from numerator and denominator
of $1-\cos(m-1)t$ which vanishes when $(m-1)t$ is an even multiple of $\pi$.  This corresponds to those values of $t$ for which the two ends of
the chord joining $(\cos t, \sin t)$ and $(\cos mt, \sin mt)$ actually {\em coincide}, but because cancellation takes place
in both numerator and denominator it is a ``removable singularity'' and the parametrization (\ref{eq:xy}) remains valid at these points. Indeed these are exactly the points
at which the envelope is {\em tangent} to the circle. To count them 
$(m-1)t=2n\pi$ ($n$ an integer) means $t=\displaystyle{\frac{2nb\pi}{|a-b|}}$ (recall
$t\ge 0$) and the values of $n$ which make this angle between 0 and $2b\pi$ are
$n=0,1,...,|a-b|$.  Thus there are $|a-b|$ points of tangency between the envelope and
the circle, proving the second part of (ii) This is also illustrated in Figure~\ref{fig:exs1}.

\subsection{Cusps}
It is clear that with $m=\frac{a}{b}$ in its lowest terms, and $b>0$, the functions giving $x$ and $y$ will begin to repeat when $t=2b\pi$.  (This is not the same as saying
that the smallest $t$ with $X(t)=X(0)$ and $Y(t)=Y(0)$ is $t=2b\pi$: the envelope may, as in (iii) of the above theorem, 
intersect itself for smaller values of $t$.) For $m$ an integer
this is of course just a period of $2\pi$.  Thus we need to consider in general the interval $0\le t < 2b\pi$ when discussing cusps or self-crossings of the envelope.
In practice we will regard parameter values $t$ as defined modulo $2b\pi$.

A cusp forms on the curve $(X(t), Y(t))$ when its velocity vector $(X'(t), Y'(t))$ is zero, the prime $'$ here representing differentiation with respect to $t$.  These require
both $\sin t+\sin mt=0$ and $\cos t + \cos mt=0$, which in particular implies
that $\sin t \cos mt=\cos t \sin mt$, that is $\sin(m-1)t = 0$ or $(m-1)t=n\pi$ for an
integer $n$. Substituting $mt = n\pi  + t$ into the equations $X'(t)=Y'(t)=0$ however
requires $\cos n\pi = -1$ in both cases, so in fact $n$ has to be {\em odd}. The even
values  correspond with the two points on the unit circle coinciding, a case
already mentioned above. 

Thus for the formation of a cusp $|m-1|t$ takes values
$(2k-1)\pi$ for $k=1, 2, \ldots$. When $m$ is an integer, we are looking for values of $0\le t <2\pi$ and the largest value
of $k$ is then $|m-1|$.  In general when $m=\frac{a}{b}$ in lowest terms we have $0\le t<2b\pi$ and the largest value
of $k$ is $|a-b|$.  Hence there are $|a-b|$ cusps in general, proving (ii) of the theorem.

\begin{rem}\label{simple}{\rm
These cusps are always what is called ``simple cusps'', which means that by a change
of coordinates in the plane close to any one cusp it can be reparametrised as
the ``normal form'' $(t^2, t^3)$.  The standard test for a simple cusp on a
curve parametrised as $(X(t),Y(t))$ is that 
$X'=Y'=0$ but $X''Y'''-X'''Y''\ne 0$ where the derivatives
are evaluated at the cusp point.
We have more to say about this condition in Lemma~\ref{lemma:swallowtail} below.
 (There is some information about reduction to
normal form in \cite[pp.154--159]{CS}; more formal treatments are found in
more advanced books such as \cite[Ch.8]{Mon}. Note that a simple
cusp is referred to as an {\em $A_2$ singularity} of a plane curve.)

In the present case $X''Y'''-X'''Y'' (m\ne 0, -1)$ is zero precisely when $m\cos(m-1)t
+m^2-m+1=0$. However when $\cos(m-1)t=-1$ this requires $m=1$, a case
already eliminated. Thus all cusps are simple.
}
\end{rem}

\section{General envelopes of lines}\label{s:general}

 Let
\begin{equation}
F(t,x,y)=A(t)x+B(t)y+C(t),
\label{eq:Fgeneral}
\end{equation}
 so that $F=0$ is a family of straight lines, parametrized by $t$,
with envelope given by $F=0, \ F_t = A'x + B'y + C'=0$, using here and in what follows subscripts to denote partial derivatives.
Solving for $x=X(t), \ y=Y(t)$ gives (omitting the variable $t$)
\[ X = \frac{BC'-B'C}{AB'-A'B}, \ \ Y= \frac{A'C-AC'}{AB'-A'B},\]
provided of course that $AB'-A'B\ne 0$.

For the envelope considered above in \S\ref{s:1circle} (see (\ref{eq:F1circle})), $A=\sin mt-\sin t, B = -\cos mt+\cos t, C = -\sin(m-1)t.$
Calculation then shows that
$AB'-A'B= (m+1)(1-\cos(m-1)t)$ which is zero when $(m-1)t$ is an even multiple of $\pi$. These are exactly the points for which
the two ends of the chord joining $(\cos t, \sin t)$ to $(\cos mt, \sin mt)$ {\em coincide}, corresponding to lines in the family which (as limits) are tangent to
the circle. Fortunately this does not invalidate the parametrization (\ref{eq:xy})
and indeed the envelope is actually smooth at these points, the cusps occurring at intermediate values of $t$ where $(m-1)t$ is an {\em odd} multiple of $\pi$.
It is worth recording the following well-known lemma.
\begin{lemma}
In the above notation, suppose that $AB'-A'B\ne 0$ and let $(X(t), Y(t))$ be the resulting parametrization of the envelope,
At an envelope point $(X(t),Y(t))$ with corresponding parameter $t$, given by solving $F(t,X(t),Y(t))=F_t(t,X(t),Y(t))~=~0$, the second derivative $F_{tt}(t,X(t),Y(t))$ vanishes   if and only if the envelope is singular, that is
 $X'(t)=Y'(t)=0$.
 \label{lemma:Ftt}
\end{lemma}
Note that $F_{tt}(t,x,y)=A''x+B''y+C''$ so $F_{tt}(t,X(t),Y(t))= A''X(t)+B''Y(t)+C''$,
a function of $t$ only, and similarly for higher derivatives of $F$.

\noindent
{\bf Proof} \ We have $F(t,X(t),Y(t))= 0$ and $F_t(t,X(t),Y(t))=0$ identically as functions of $t$. Differentiating
these with respect to $t$ we get $F_t+F_xX'+F_yY'=0=F_{tt}+F_{tx}X'+F_{ty}Y'=0$.  It is clear that if $X'=Y'=0$ then $F_{tt}=0$.   For the converse,
assume $F_{tt}=0$ at an envelope point, so that $F_t=0$ too.  We have $F_x=A, F_y=B, F_{tx}=A', F_{ty}=B'$ so that the two equations for $X',Y'$ are
$AX'+BY'=0, A'X'+B'Y'=0$ and these imply $X'=Y'=0$ since $AB'-A'B\ne 0$.  \hfill$\Box$

\smallskip\noindent
Thus away from points where $AB'-A'B=0$ we can detect singular points of the envelope by the equation $F_{tt}=0$.  When $AB'-A'B=0$ then in general
we might expect the envelope to ``go to infinity'' since the denominator of $X$ and $Y$ vanishes.  For the envelope studied in
 \S\ref{s:1circle} numerator and denominator always vanish together and these ``singularities'' are removable. The envelope {\em never}
 goes to infinity and its genuine singular points are cusps as determined in Theorem~\ref{th}.
 In the next section we shall meet an example where the equal denominators of $X$ and $Y$ can vanish with neither numerator, or one numerator, vanishing,
though never both at once. When one numerator vanishes this means that the envelope ``goes to infinity'' in the direction of the $x$- or
 the $y$-axis, and when neither vanishes it means that the envelope goes to infinity in some other direction.
 
 \medskip
 
 Naturally singularities of envelopes $(X(t),Y(t))$, detected by the conditions $X'=Y'=0$, are not always simple cusps, in the sense that they do not always ``look like'' the
 cusp $y^2=x^3$. In the case of
 the singularities present in one-circle envelopes they are in fact simple cusps: a suitable smooth change of coordinates in the plane near to the cusp
 will transform it into the standard form $(t^2, t^3)$. The  test for this is
 that, writing $(X(t),Y(t))$ for a parametrization of the envelope, such as (\ref{eq:xy}),
 we have $X'(t_0)=Y'(t_0)=0$, but $ X''(t_0)Y'''(t_0)\ne X'''(t_0)Y''(t_0)$, where $t_0$ gives the position of the cusp. 
 The latter condition is equivalent in fact to $F_{ttt}(t,X(t),Y(t)) \ne 0$
 at $t=t_0$.
Calculation shows that either of these conditions reduces to $2\cos((m-1)t_0)+1=0$,
which is incompatible with the condition $\sin((m-1)t_0)=0 $ for $F_{tt}=0$, so
all cusps are simple. 
 
 In the next section we shall meet higher singularities. Simple cusps are
 {\em stable} in the sense that small changes of continuously varying parameters
 may move their position but will not destroy the cusp itself.  For the one-circle
 envelopes there are no such parameters to vary ($m$ is restricted to being rational) but 
 we shall now introduce one, and later two, parameters of this kind to vary.

\section{An envelope from two circles}\label{s:2circles}

There is an interesting extension of the above construction to produce a {\em 
one-parameter family of envelopes} in which changing the
parameter which increases
or decreases the number of cusps and self-intersections. However, as we shall see,
the construction is too symmetric to exhibit the most general behavior, and we 
need to break the symmetry by means of a 2-parameter family to reveal the
full structure of the singular points which arise.

We take two circles centered at the origin, of radii 1 and $r >1$, a fixed
rational number $m$, and join the point
$(\cos t, \sin t)$ of the first circle to the point $(r\cos mt, r\sin mt)$ of the second
with a straight line. This line has equation $F(t,x,y,r)=0$, where
\begin{equation}
 F(t,x,y,r) = x(r\sin mt-\sin t)- y(r\cos mt -\cos t) - r\sin (m-1)t.
\label{eq:F-2circles}
\end{equation}
Here $r$ is a continuous parameter, given a definite value 
for each example, and $t$ tells us
which line we are considering.

We shall later describe $F$ as a {\em three-parameter family of functions of the
variable $t$}. Note that $m$ is restricted to rational values, so for example cannot 
be differentiatied, and is regarded as fixed.

The envelope of the lines $F=0$, obtained from $F=\partial F/\partial t = 0$, is
(suppressing $r$ from the notation for $X,Y$ now)
\begin{eqnarray}
X(t)&=& \frac{r\left[\cos mt(r\cos(m-1)t -1) + m\cos t(\cos(m-1)t - r)\right]}{r(m+1)\cos(m-1)t - (mr^2+1)},\nonumber \\
&& \label{eq:xy-2circles} \\
Y(t)&=& \frac{r\left[\sin mt(r\cos(m-1)t-1) + m\sin t(\cos(m-1)t - r)\right]}{r(m+1)\cos(m-1)t - (mr^2+1)}. \nonumber
\end{eqnarray}
Note that in this case the factor $\cos((m-1)t)-1$, which was cancelled from numerator and denominator when deducing the
parametrization (\ref{eq:xy}) in the one-circle case, appears here explicitly: when $r=1$ the
denominator in (\ref{eq:xy-2circles}) becomes exactly $(m+1)(\cos (m-1)t -1)$.

We shall assume, as before, that $m=\frac{a}{b}$, a rational number in its lowest terms, with $b>0$. We shall exclude certain values of $m$, namely $-1,0,1$, and sometimes
$2$ and $\ha$. 

\begin{rem}\label{rem:T}
{\rm
The two-circle envelope above can also be constructed from the family of straight lines joining the point $(\cos bT, \ \sin bT)$ on the unit circle
to the point $(r\cos aT, \ r\sin aT)$ on the concentric circle of radius $r$: just write $t=bT$ so that $mt=(a/b)bT=aT$. The range of values of $T$
is $0\le T <2\pi$. The formula corresponding to~(\ref{eq:F-2circles}) is
\[ \widetilde{F}(x,y)=x(r\sin aT - \sin bT)-y(r\cos aT - \cos bT) - r\sin(a-b)T. \]
}
\end{rem}

\begin{figure}[!ht]
\begin{center}
\includegraphics[width=6.8in]{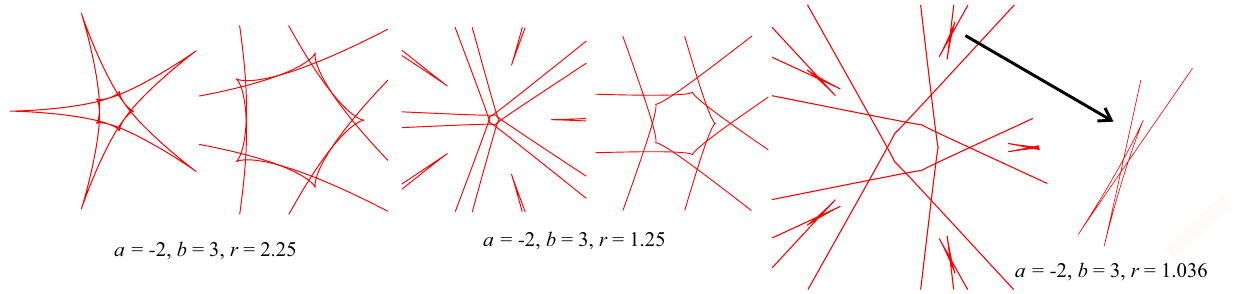}
\end{center}
\vspace*{-0.25in}
\caption{\small Left pair:  $r>\frac{1}{|m|}= 1.5$ and the envelope is finite;  the central part is enlarged.
Centre pair: the envelope goes to infinity 10 ($=2(b-a)$) times, because $r<\frac{1}{|m|}$ (see Proposition~\ref{prop:env-infty}); the central part of the
envelope is again enlarged for clarity. Right pair: $r<\left|\frac{m-2}{2m-1}\right| = \frac{8}{7}$ which results in 10 additional cusps, as in Proposition~\ref{prop:cusps-2circles}.
One of the configurations of three cusps is enlarged at the far right. }
\label{fig:env-infty}
\end{figure}

Note that the denominator in (\ref{eq:xy-2circles}) can be zero, indicating that the envelope has ``gone to infinity''.  This occurs when
\begin{equation}
\cos(m-1)t = \frac{mr^2+1}{r(m+1)}.
\label{eq:env-infty}
\end{equation}
The denominator in {\em this} expression is not zero! Values of $t$ exist satisfying (\ref{eq:env-infty}) if and only if the right-hand side lies in the closed interval
$[-1,1]$.
Some work with inequalities, using $r>1$ and $m \ne -1, 0, 1$, shows that this is equivalent to $1<r\le \frac{1}{|m|}$
and in particular is only possible if $|m|<1$. Summing this up:
\begin{prop}
The 2-circle envelope parametrized by {\rm (\ref{eq:xy-2circles})}, where $r > 1, \ m=\frac{a}{b}, b>0$ is in lowest terms,
 goes to infinity if and only if $|m|<1$ and $1<r\le\frac{1}{|m|}$. The values of $t$ for which
this happens are given by {\rm (\ref{eq:env-infty})} which has $2|b-a|$ solutions for $0\le t < 2b\pi$, except for $r=\frac{1}{|m|}$ when the number is $|b-a|$. 
\hfill$\Box$
\label{prop:env-infty}
\end{prop}
There is an illustration of this in Figure~\ref{fig:env-infty}.

\subsection{Cusps on the 2-circles envelope}\label{ss:cusps2circle}
\medskip\noindent
When counting the cusps on the 2-circles envelope it is helpful to use Lemma~\ref{lemma:Ftt}.
At points given by (\ref{eq:xy-2circles}), we have
\begin{equation}
F_{tt} = rm\sin(m-1)t\left( m-2+r(m+1)\cos(m-1)t+r^2(1-2m)\right)
\label{eq:cusps-2circles}
\end{equation}
and cusps (or possibly higher singularities: see \S\ref{s:higher}) occur when this is zero. 

The factor $\sin(m-1)t$ is equal to zero if $|m-1|t$ takes values $k\pi$ 
for $k=0, 1, \ldots$. Thus, taking $m=\frac{a}{b}, b>0$ in lowest terms and looking at the range $0\le t<2b\pi$, there are $2|a-b|$ 
solutions for this term, corresponding to $2|a-b|$ ``general cusps'' for any $r>1$.

\noindent The additional factor of (\ref{eq:cusps-2circles}) equals zero when
\begin{equation}
\cos(m-1)t = \frac{(2m-1)r^2 + (2-m))}{r(m+1)}=M \mbox{ say.}
\label{eq:cusps-extra}
\end{equation} 
Values of $t$ satisfying this expression exist if and only if $M$ lies within the interval $[-1, 1]$. Evaluating the 
resulting inequalities gives the result that extra cusps occur in the range $1<r\le |\frac{m-2}{2m-1}|$ if $|m|<1$. Conversely, for $|m|>1$, 
there are no extra cusps besides the $2|a-b|$ general cusps for all $r>1$. When $r$ takes the extreme value $|\frac{m-2}{2m-1}|$ then
it is easy to check that $\cos(m-1)t=\pm1$, so that $\sin(m-1)t=0$ and the values of $t$ have coincided with those for the ``general cusps'' above, possibly changing the
nature of these cusps.

\begin{rems}{\rm
(1) \ The values of $t$ in (\ref{eq:cusps-extra}) cannot coincide with those in (\ref{eq:env-infty}) (given as usual $m\ne 1, r> 1$), but the cusps given by
$\sin (m-1)t=0$ can be at infinity, in fact when $|m|=1/r$.  An example occurs at the intermediate value $r=1.5$, between the illustrations in Figure~\ref{fig:env-infty}, left
and center,
when all these cusps (the outer five in the left-hand diagram) have ``gone to infinity''.

\sn
(2) \ Clearly $m=\frac{1}{2}$ is slightly special in the above discussion: with $a=1, b=2$ there are $2|a-b|=2$ additional cusps for all values of $r>1$.  The envelope
goes to infinity in this case just for $r\ge 1/|m| = 2$, according to Proposition~\ref{prop:env-infty}, and when $r=2$ one of the cusps goes to infinity.
}
\end{rems}

\begin{figure}[!ht]
\begin{center}
\includegraphics[width=2in]{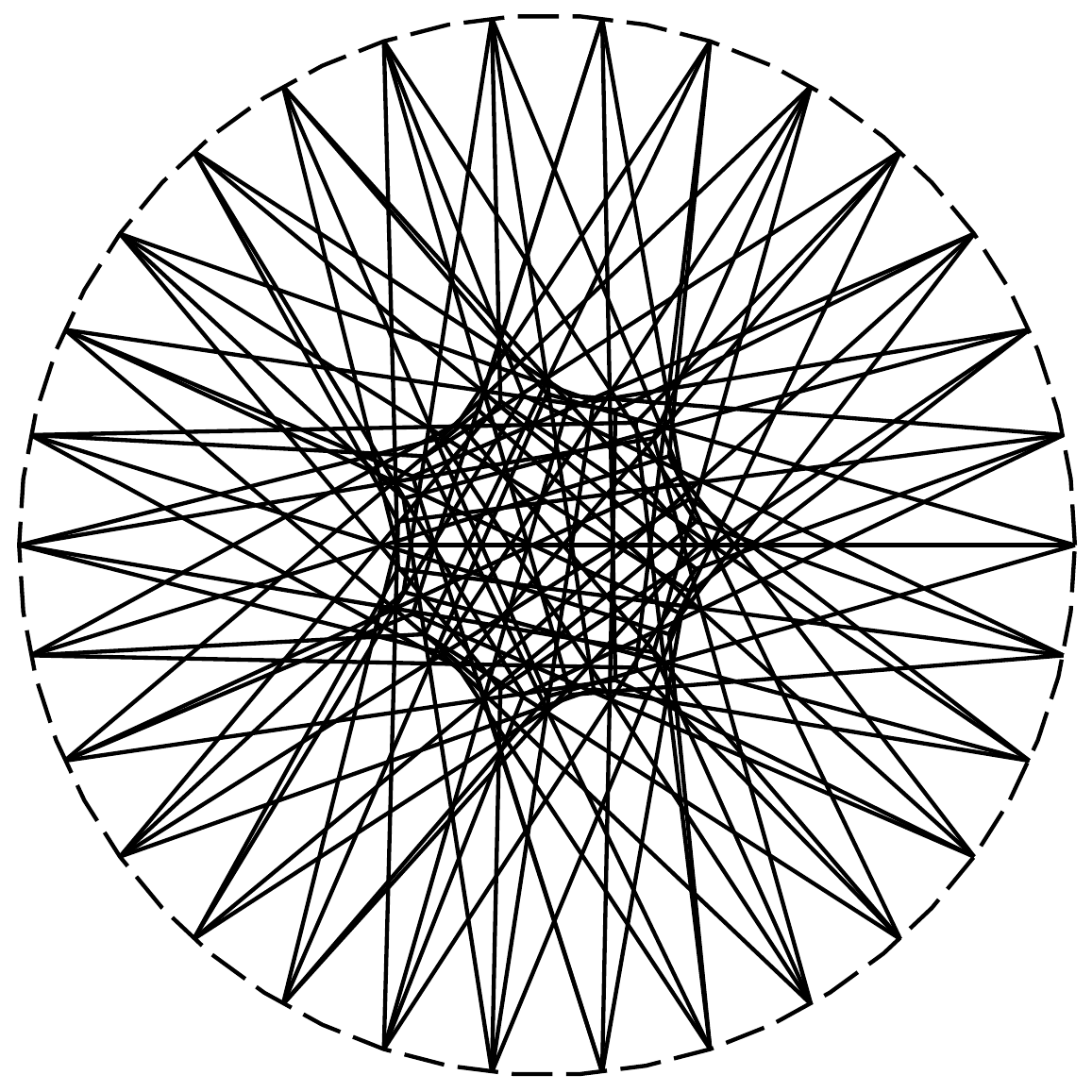}
\hspace*{0.31cm}
\includegraphics[width=2in]{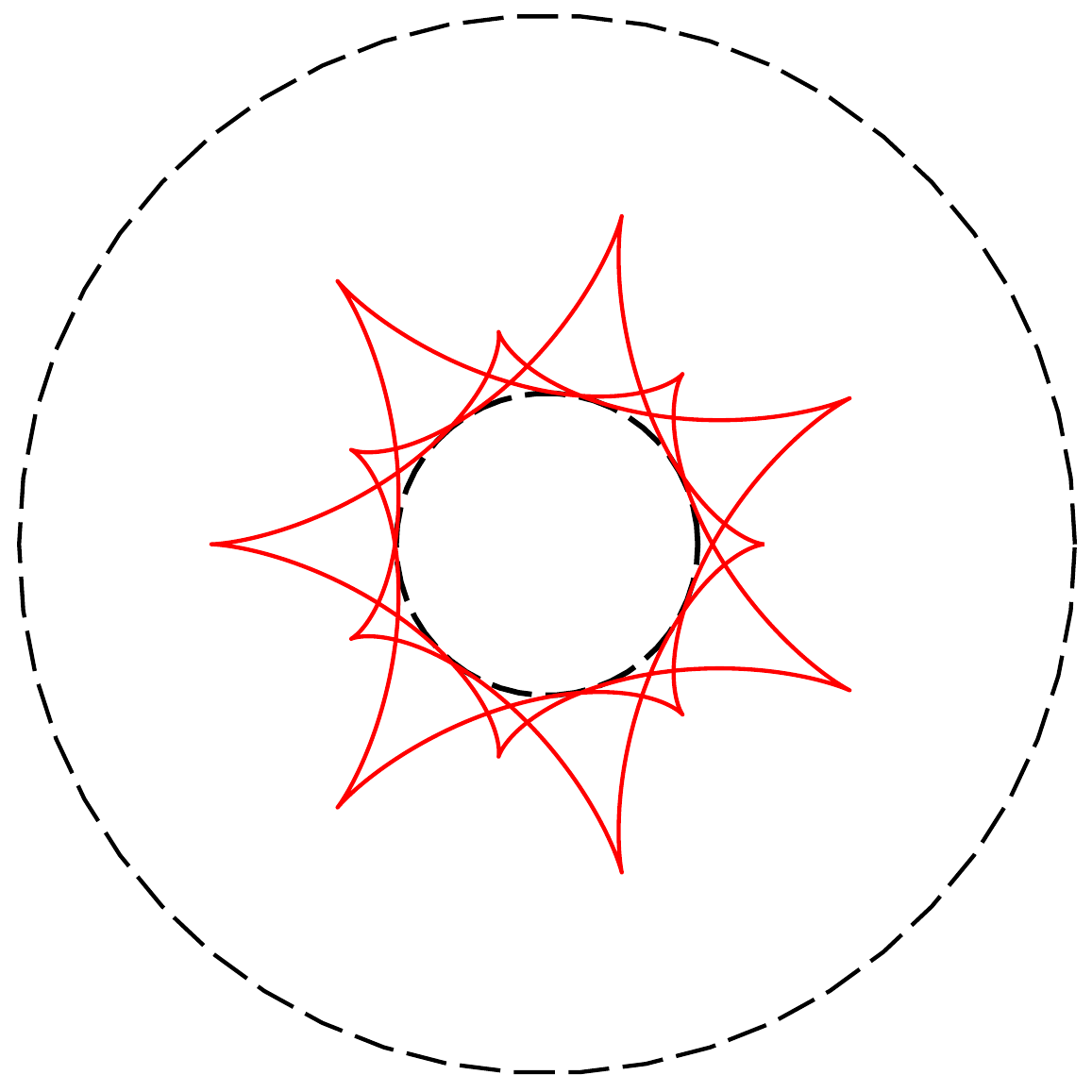}
\hspace*{0.31cm}
\includegraphics[width=2in]{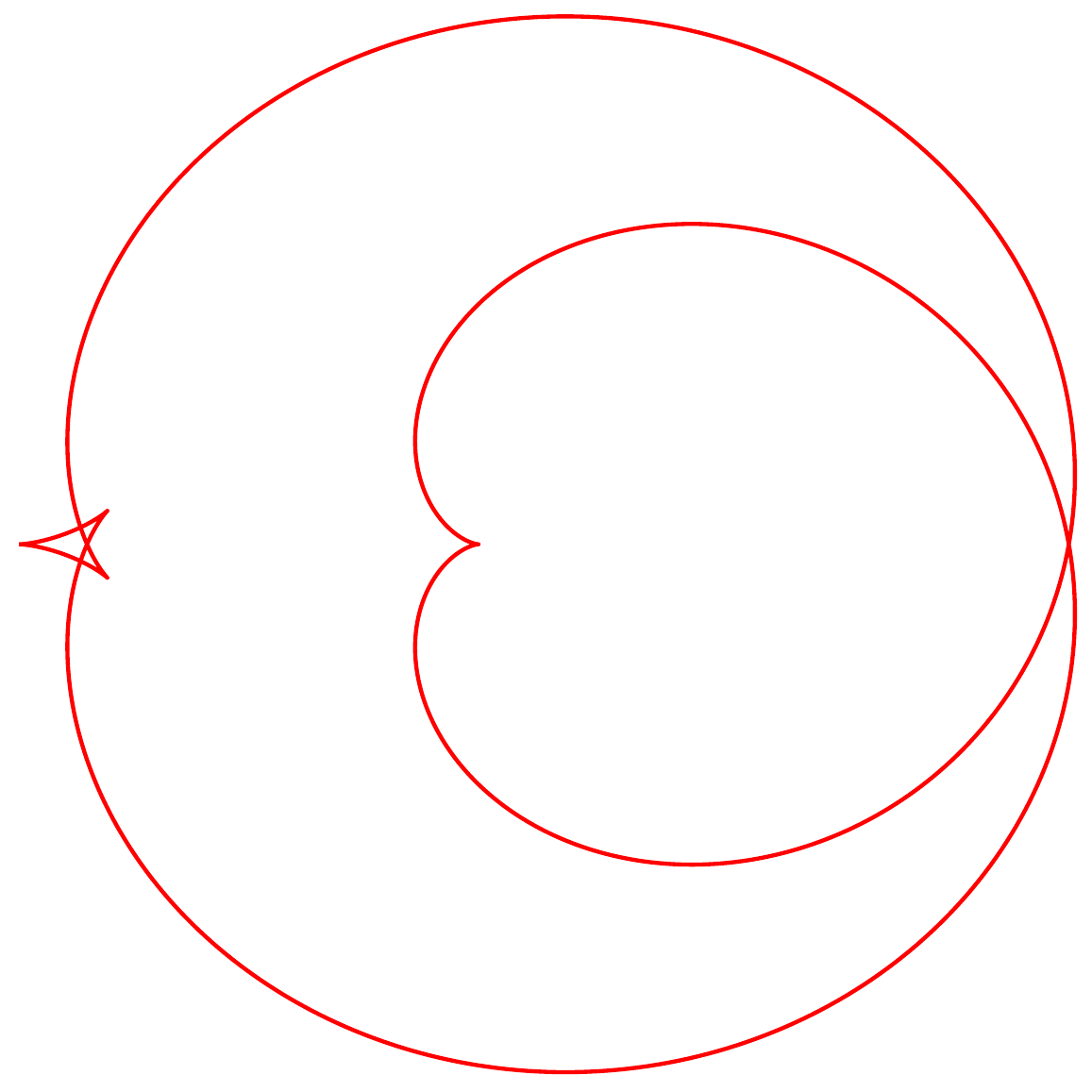}
\end{center}
\caption{\small Left and center: the case $a=-4, b=3, r=3.5$, with $2|a-b|=14$ cusps. On the left the lines forming the envelope
are drawn and also the circle radius $r$ is drawn dashed. In the center the envelope is the solid curve and the
circles of radii 1 and $r$ are drawn dashed.  Right: the case $a=3, b=4, r=2$ with 4 cusps, which is  $4|a-b|$, as in (iia) of Proposition~\ref{prop:cusps-2circles}.
The collection of three nearby cusps at the left of this envelope is called a {\em butterfly configuration}. As $r\to |2-m|/|2m-1|=2.5$ this configuration collapses to a single 
point which is a {\em butterfly singularity}.  This is not a simple cusp and is discussed further in \S\ref{s:higher}. At $r=2.5$ there are now $2|a-b|=2$ cusps as in 
(iib) of Proposition~\ref{prop:cusps-2circles}.
See also Figure~\ref{fig:butterfly1}.}
\label{fig:2circles1}
\end{figure}

\noindent
From the above discussion we have the following.
\begin{prop}
The two-circle envelope parametrized by {\rm (\ref{eq:xy-2circles})}, where as usual $m=\frac{a}{b}$ is in its lowest terms, $r>1$ and the range of $t$ is $[0, 2b\pi)$, has cusps as follows.

\sn
{\rm (i)} \ $2|a-b|$ cusps always (we call these the general cusps); \\
{\rm (iia)} \ when $|m|<1$ and  $1<r<\left|\frac{2-m}{2m-1}\right|$,  $2|a-b|$ cusps in addition to (i) {\rm [}two additional cusps for all $r>1$ when $m=\frac{1}{2}${\rm ]}
See figure~\ref{fig:2circles1}, right; 
\\
{\rm (iib)} \ when $|m|<1$ and $r=\left|\frac{2-m}{2m-1}\right|$,  $2|a-b|$ cusps 
altogether.  Not all these cusps may be simple; this is explored in \S\ref{s:higher}.
\hfill$\Box$
\label{prop:cusps-2circles}
\end{prop}

\begin{figure}[!ht]
\begin{center}
\includegraphics[width=1.5in]{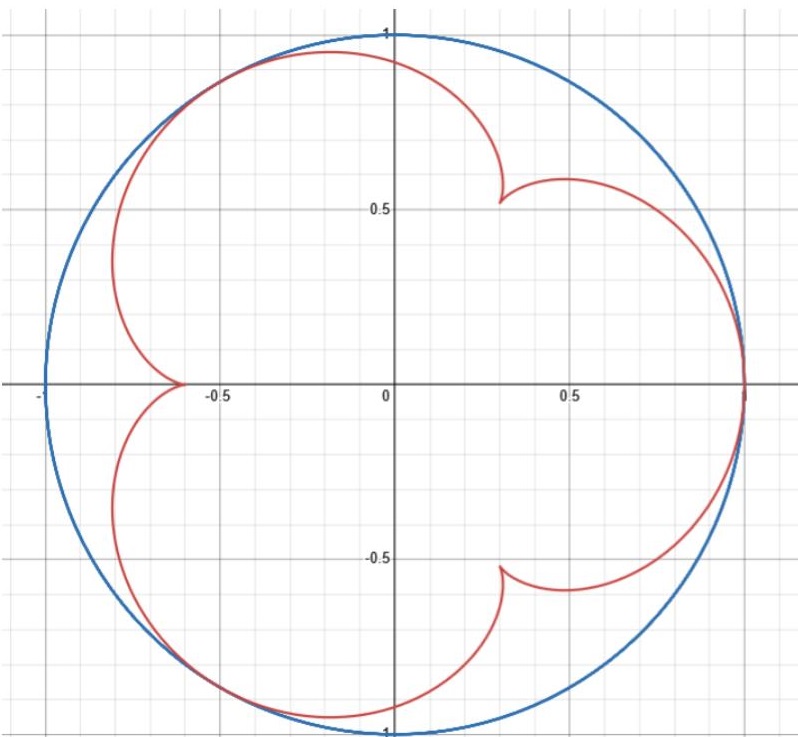}
\includegraphics[width=1.5in]{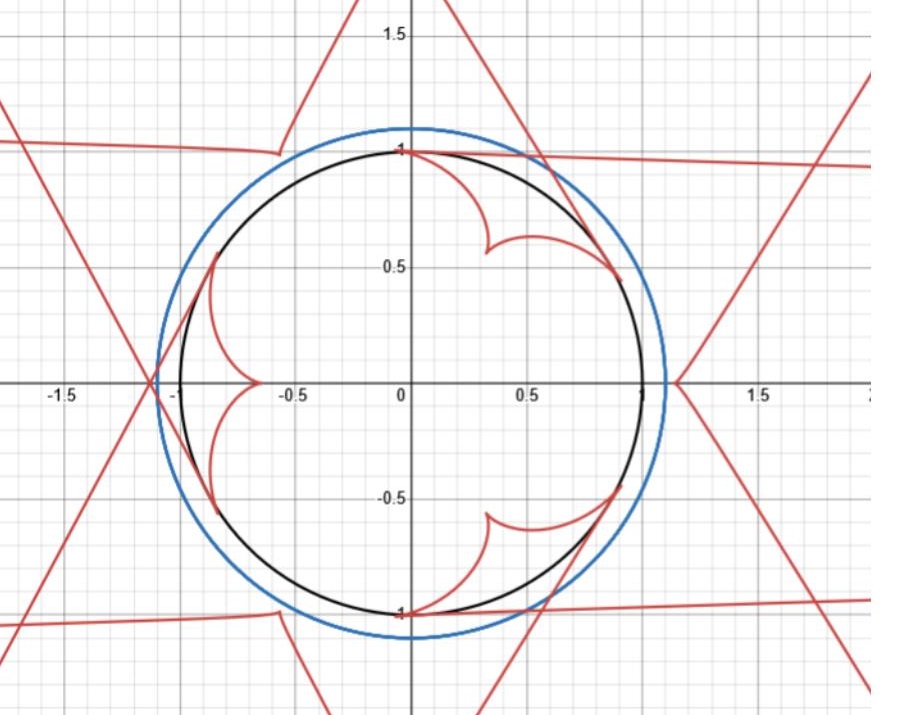}
\includegraphics[width=1.5in]{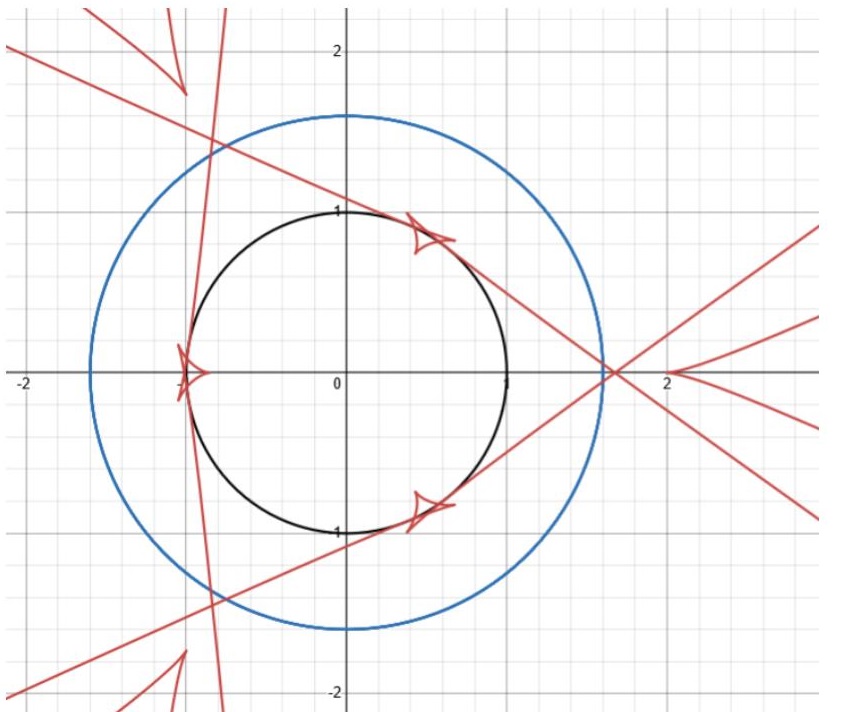}
\includegraphics[width=1.5in]{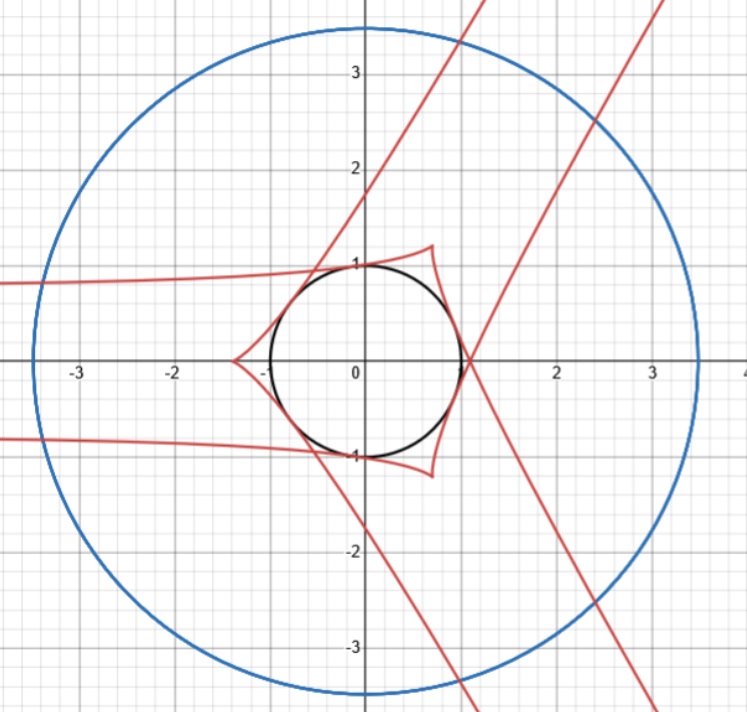}
\end{center}
\vspace{-0.3in}
\caption{\small A discontinuous change for $m=1/4$ from $r=1$ to $r=1.1$ and then on 
to $r=1.6$  and $r=3.5$. The last transition shows three butterfly configurations which
collapse to three (non-simple) cusps. To reveal the full structure of these cusps we need
to introduce a second parameter as in \S\ref{s:higher} and Figure~\ref{fig:butterfly1}.}
\label{fig:2circles-m=1}
\end{figure}

Some further examples are given in Figures~\ref{fig:2circles1} and \ref{fig:2circles-m=1}.

\begin{figure}
\begin{center}
\includegraphics[width=3.8in]{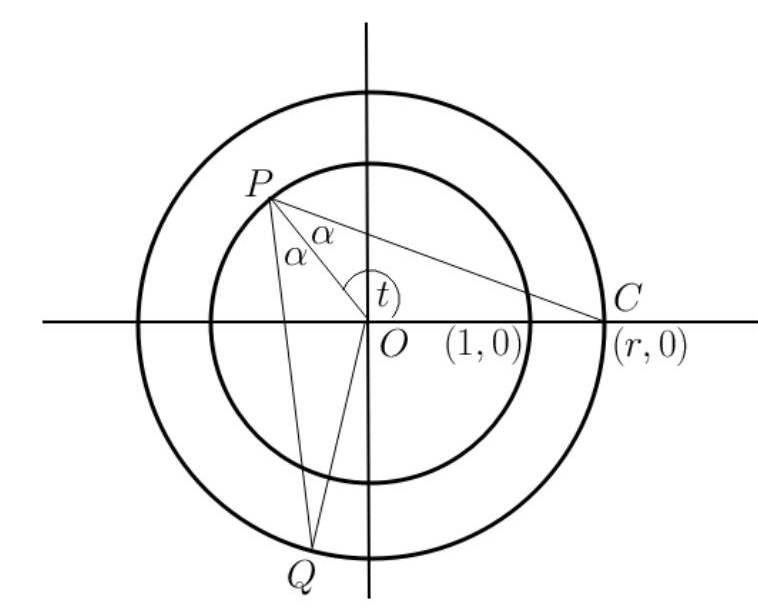}
\includegraphics[width=3in]{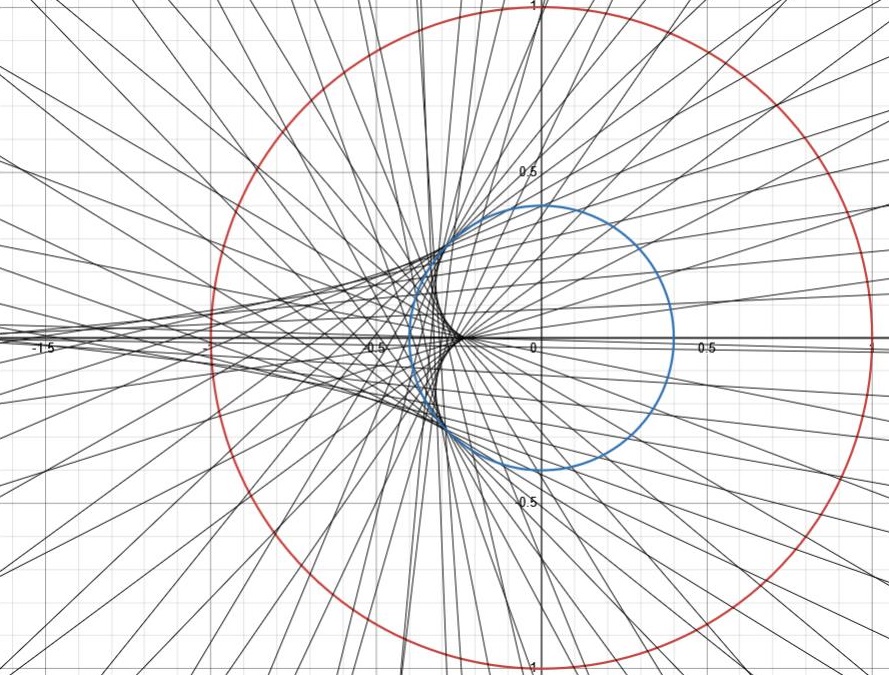}
\vspace{-0.2in}
\end{center}
\caption{\small Left: identification of a caustic by reflexion in the unit circle with source
at $(r,0)$ with the 2-circle envelope taking $m$ to be 2. Right: envelope of reflected rays with $r=0.4$; thus the source of light is at the point $(0.4,0)$ on the inner circle and the
full mathematical caustic is drawn, irrespective of whether the reflected
rays focus {\em inside}
the mirror (the outer circle).}
\label{fig:caustic}
\end{figure}

\begin{rem} The special case $m=2$ (or $m=\ha$).\label{rem:caustic}
{\rm
The 2-circle envelope of lines, in this special case, has an interpretation as a
different envelope of lines, in fact as the envelope of rays reflected from a circular
mirror of radius 1 with a point light source at  $(r,0)$ where $r>0$ is  the
radius of the second circle. (By symmetry it is enough to consider a light
source on the positive $x$-axis.) This is the {\em caustic by reflexion} of the unit
circle with a finite light source. For further information on caustics see for example
\cite{BGG,LB}; there are also articles and videos on the internet under headings
such as `caustic' (in the mathematical, not chemical, sense).

Figure~\ref{fig:caustic}, left, shows a light ray $CP$
from a point source at $(r,0)$ reflected at the point $P$ on the unit circle 
(the ``mirror'') and striking
the concentric circle of radius $r$ at $Q$. Note that if the angle $t$ is small 
in magnitude then
the light ray will in practice bounce off the unit circle and not focus at all: we have to
produce it ``backwards'' inside the unit circle for focussing to take place. Thus we are
really considering the ``mathematical caustic'' which is a completion of the ``physical
caustic''. In practice light rays may be reflected many times inside the mirror;
many studies have considered these multiple caustics, for example \cite{BT}.

The triangles $POC$ and $POQ$ are 
congruent so that the point $Q$ has coordinates $(r\cos 2t, \ r\sin 2t)$, that is
$PQ$ is the same line as in the construction of the 2-circle  envelope with $m=2$.
Any lingering doubts that this always works, even  if $r<1$,
 can be answered by finding the
equation of the reflected ray:
\[  \sin t(2r\cos t -1)x + (\cos t -r\cos 2t)y - r\sin t = 0  \]
which always contains the point $(r\cos 2t, \ r\sin 2t)$ using the standard double angle
formulas.
}
\end{rem}

\section{Swallowtails and Butterflies} \label{s:higher}
The cusps on an envelope constructed from two concentric circles 
call for some further analysis. One such envelope is shown in the central diagram of
Figure~\ref{fig:butterfly1} (the circles themselves are not shown), where 
$m=\frac{3}{4}$ and  $r=\frac{5}{2}$.
 In terms of Proposition~\ref{prop:cusps-2circles}(iib)
this is an extreme value of $r$. The left and right cusps occur for $t=0$ and
$t=4\pi$ respectively and the initial terms of the Taylor expansions of the 
envelope curves $(X(t),Y(t))$ close to the two cusps in this example are as follows:
\begin{equation}\label{eq:taylor}
 \mbox{Left cusp: } \left(-\frac{5}{7} - \frac{25}{1792}t^4, \ \ -\frac{5}{768} t^5\right); \ \mbox{ Right cusp: } \left( -\frac{5}{23} - \frac{75}{1058}(t - 4\pi)^2, \ \ 
-\frac{25}{644}(t - 4\pi)^3)\right).
\end{equation}
The right cusp is therefore a simple cusp, reducible to $(t^2, t^3)$ by a change of
coordinates, but the left cusp is not: it is
 a ``$(4,5)$'' cusp which is called a {\em butterfly singularity}. Figure~\ref{fig:2circles1},
 right, shows the effect on this singularity of changing $r$ slightly: the butterfly
 opens its wings into a 3-cusped curve (all these cusps are in fact simple),
 while the simple cusp remains simple.

In what follows we shall use higher derivatives of the defining equation 
 (\ref{eq:F-2circles}) of the 2-circle envelope to distinguish the various cusps.
 (These conditions apply to any envelope of lines, not just the 2-circle envelope.)
 Recall that the envelope itself is determined by $F=F_t=0$. The condition
 for a singular envelope is that in addition $F_{tt}=0$ (evaluated as usual at
 $(t,X(t),Y(t))$), which from Lemma~\ref{lemma:Ftt} is equivalent to $X'=Y'=0$. 
 First we introduce a standard example of something simpler 
 (but not simple), unfolding  a  ``(3,4) cusp'' or swallowtail point.
 
 \subsection{A standard swallowtail}\label{ss:standard-swallowtail}
\begin{figure}[h!]
\begin{center}
\includegraphics[width=2.2in]{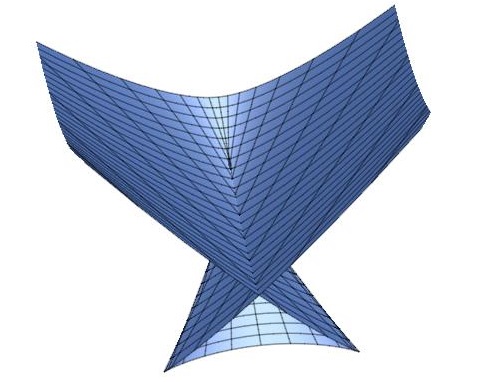}
\includegraphics[width=4.1in]{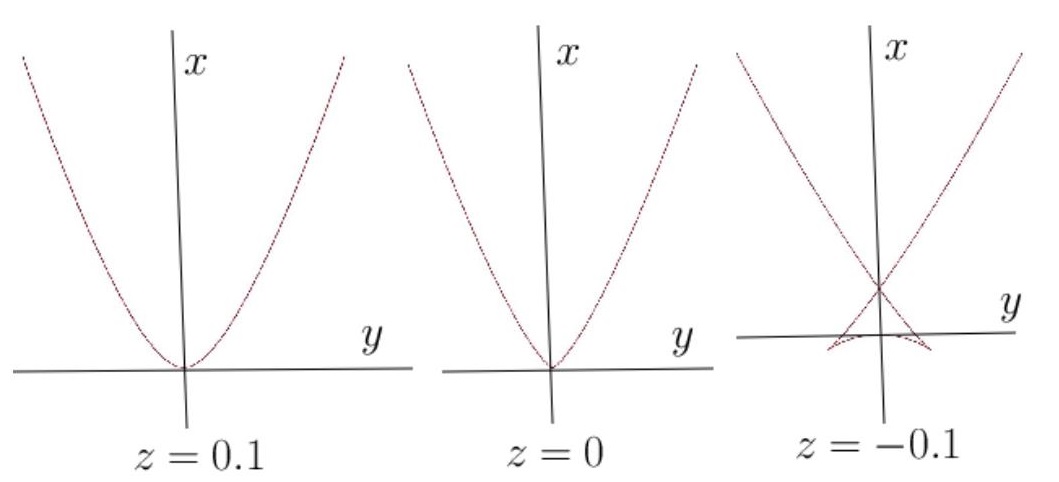}
\end{center}
\caption{\small Left: a swallowtail surface in $(x,y,z)$-space, 
obtained as the discriminant of
the family $G(t,x,y,z)=t^4+x+yt+zt^2$. Right: slices 
$z=$ constant of this surface from ``back'' to ``front''. The reversal of the usual
axes arises from using the natural progression $x,y,z$ in the definition of $G$.
The sequence displays a ``swallowtail transition'' but a similar transition, where two
cusps and a crossing merge and evolve into a smooth curve, does not occur on 2-circle envelopes without
introducing a second parameter besides $r$  (see Proposition~\ref{prop:higher}
and the comments following it). The ``swallowtail
singularity'' itself occurs in the center diagram. It looks like a kink in the curve, a point where the first derivative exists as a limit but the second does not. 
 }
\label{fig:swallowtail}
\end{figure}

 A basic example of  a 3-parameter family of functions of $t$,
  not arising from 2-circle envelopes, is
\begin{equation}\label{eq:G}
G(t,x,y,z)=t^4+x+yt+zt^2
\end{equation}
for which the discriminant set $G=G_t=0$
is parametrised 
$x(t)=3t^4+zt^2, \ y(t)=-4t^3-2zt$. (Thus $z$ here plays an analogous role
to $r$ in the 2-circle envelope case.) This set is illustrated in $(x,y,z)$
 space in Figure~\ref{fig:swallowtail}, left, 
and the slices $z=$ constant near to $z=0$ in Figure~\ref{fig:swallowtail}, right.
When $z=0$ the parametrisation of the curve is $x=3t^4,\ y=-4t^3$, a ``(3,4)''
singularity or swallowtail singularity.
 As $z$ goes through 0, the swallowtail singularity unfolds to
reveal two ordinary cusps and a crossing in one direction $z<0$ and a smooth
curve in the other direction $z>0$.

Note that in this example, $G_{tttt}(0,0,0,0)=120 \ne 0$. We shall see below
that this property distinguishes $(3,4)$ singularities on  2-circle envelopes
from the standard example above. It is what prevents a swallowtail point
on a 2-circle envelope from being fully unfolded by changing the radius $r$.
We shall see in \S\ref{ss:unf} how adding an extra parameter to $F$, thereby
breaking the symmetry of the construction, both swallowtail and higher
``butterfly'' singularities exhibit full unfoldings.

\subsection{Further details of the 2-circle construction}

 
 In general terms, the more derivatives 
 of $F$ with respect to $t$ which vanish, the ``higher'' the cusp. 
 A simple cusp has $F=F_t=F_{tt}=0$ but $F_{ttt}\ne 0$. Reference****
 Here we shall check the next case.
 \begin{lemma}\label{lemma:swallowtail}
 Suppose that $F=F_t=F_{tt}=0$ at $(t,X(t),Y(t))$ so that, as above,
 the envelope of lines $F=0$ is singular at this point.  Then $F_{ttt}=0$ too if and
 only if $X''Y'''-X'''Y''=0$, which is the condition for a ``higher cusp'', typically
 a ``swallowtail point''.
 \end{lemma}
 {\bf Proof} \ The equations $F(t,X(t),Y(t))\equiv 0$,
 that is $aX+bY+c\equiv 0$  and $F_t(t,X(t),Y(t))\equiv 0$, that is $a'X+b'Y+c'\equiv 0$, are
 {\em identities} and can be differentiated with respect to $t$ any number of times,
 still giving zero. Differentiating the first and using the second gives $aX'+bY'\equiv 0$, which
 simply says that the line $F=0$ is tangent to the envelope curve. Differentiating this
 and using $X'=Y'=0$ gives $aX''+bY''=0$ at singular points.

 Next, differentiating $F=0$ three times, and $F_t=0$ twice, using $X'=Y'=0$ gives
 the equations
 \[ F_{ttt}+3a'X''+3b'Y''+aX'''+bY'''=0, \ \ F_{ttt}+a'X''+b'Y''=0\]
 respectively, from which we deduce $aX'''+bY'''=2F_{ttt}$ at points where $X'=Y'=0$.   
 Combining this with $aX''+bY''=0$ and using the fact that $a$ and
 $b$ cannot vanish together, gives the result of the lemma. \hfill$\Box$

At such a swallowtail point, it can be shown that by a local change of coordinates
the envelope assumes the ``normal form'' $(t^3, t^4)$. Examples of swallowtail points
are given below, though in fact the 2-circles envelope does not actually
exhibit them: we need to break some of the symmetry to observe such singularities.
(See (iib) in the Proposition below.)
The fourth derivative of $F$, that is $F_{tttt}$, can vanish in addition to the second
and third, and in that case the singularity is called a ``butterfly''. We have already met
this in (\ref{eq:taylor}) above.
(See also~\cite[p.129]{CS}, or, for the more advanced viewpoint of
``elementary catastrophes'', \cite[\S7.5]{Mon}, \cite[\S13]{BL} or \cite{PS}. There
is also much information, and good illustrations, on the internet, such as ~\cite{W}.)  

\medskip

Recall from Proposition~\ref{ss:cusps2circle} that 	cusps occur on the 2-circles
envelope when $F_{tt}(t,X(t),Y(t))=0$ where $F=0$ is the equation of the line as
in (\ref{eq:F-2circles}) and $F_{tt}$ is written out in (\ref{eq:cusps-2circles}).
Furthermore $F_{tt}=0$ if and only if $\sin(m-1)t=0$ or 
$\cos(m-1)t =M=((2m-1)r^2 + 2-m)/r(m+1)$, which is only possible when
$|m|<1$ and $1<r\le |m-2|/|2m-1|$.

Routine  but sometimes messy calculations then reveal the following details. 

\begin{prop}\label{prop:higher}
{\rm (i)} \ Suppose that $\sin(m-1)t=0$, that is $(m-1)t$ is an integer multiple of $\pi$.
Then the non-simple cusp condition $F_{ttt}(t,X(t),Y(t))=0$ (or $X''Y'''=X'''Y''$)
requires $r$ to have an extreme value $|m-2|/|2m-1|$.
In fact 
$r = (2-m)/(2m-1)$ when $\cos(m-1)t = +1$, 
and $r = (m-2)/(2m-1)$ when $\cos(m-1)t=-1$.  \\
{\rm (ii)} When $\sin(m-1)t=0$ and $F_{ttt}=0$ as in {\rm (i)}, then also
$F_{tttt}=0$, which raises the singular type of the cusp to the form $(t^4,t^5)$
as in the central figure of Figure~\ref{fig:butterfly1}. Mercifully, the next derivative
is not zero here!\\
{\rm (iii)} When $\cos(m-1)t=M$ then the
non-simple cusp condition $F_{ttt}=0$ requires that $r$ is an extreme value $|(m-2)/(2m-1)|$, but in that
case $M=\pm 1$ so $\cos(m-1)t = \pm 1$ and we are back in case {\rm (i)}, that is
$\sin(m-1)t=0$. \hfill$\Box$
\end{prop}

Note that  the
two cusps in the central figure of Figure~\ref{fig:butterfly1} are distiguished
by (i) and (ii) of the Proposition. For the left hand
non-simple cusp,
$t=0$ and $m=3/4, \cos(m-1)t=+1, r=5/2=(2-m)/(2m-1)$ so $F_{ttt}=0$,  
while for the right hand
simple cusp, $m=3/4, t=4\pi, \cos(m-1)t=-1$ but $r\ne (m-2)/(2m-1)$ so $F_{ttt}\ne 0$.

Note also that (ii) and (iii) of the Proposition show that swallowtail singularities
given by $F_{tt}=0, F_{ttt}=0, F_{tttt}\ne 0$ do not occur on 2-circle envelopes. We shall
see below that they occur when we add another parameter besides $r$ to the family
of lines defining the envelope.

We can treat the left hand side of equation (\ref{eq:F-2circles}) 
as $F(t,x,y,r)$, a 3-parameter family of functions of $t$.  For fixed $r$ the {\em discriminant set} is the envelope
in the $(x,y)$ plane as above, given by $F=F_t=0$.

We say something about the underlying reason why $G$ displays a swallowtail transition and $F$
does not, that is why the fourth derivative is important,
 in \S\ref{ss:tech}.

\subsection{A standard butterfly}\label{ss:H}
 
\begin{figure}[h!]
\begin{center}
\includegraphics[width=5in]{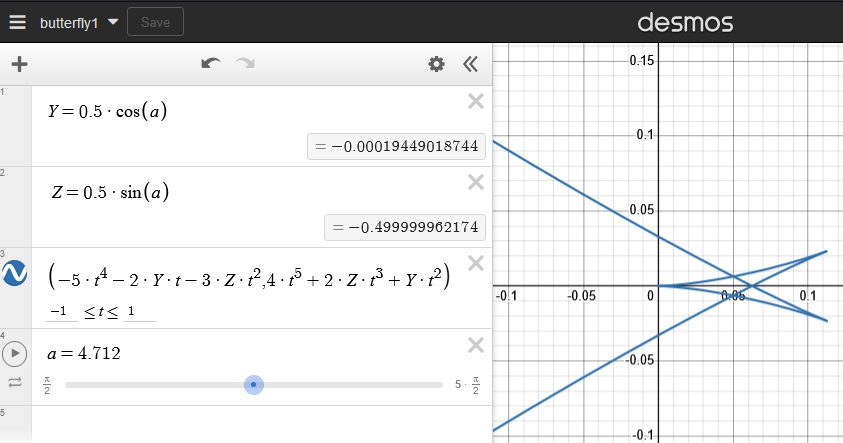}
\end{center}
\caption{\small One of the curves given by fixing $(y,z)$ (appearing here
as $(Y,Z)$) close to $(0,0)$
in (\ref{eq:butterfly-param}). This Desmos figure sets $Y=0.5\cos a, \ Z=0.5\sin a,
\pi/2 \le a \le 5\pi/2$ so that a tour of the origin $(Y,Z)=(0,0)$ occurs by allowing
the angle $a$ to travel through this interval. Clicking on the arrow against $a$ animates
this tour. The figure initially displays $a=3\pi/2\approx 4.712$.
 }
\label{fig:desmos-butterfly}
\end{figure}

 By analogy with the family $G$ in (\ref{eq:G}) for a swallowtail consider the
 family
 \begin{equation}\label{eq:H}
 H(t,w,x,y,z)=t^5+w+xt+yt^2+zt^3
 \end{equation}
 for which the discriminant $H=H_t=0$ in 4-dimensional $(w,x,y,z)$ space has parametrisation
 \begin{equation}\label{eq:butterfly-param}
 (w,x,y,z)=(4t^5+2zt^3+yt^2, \ -5t^4-2yt-3zt^2, \ y, \ z). 
 \end{equation}
This is a 3-dimensional set in 4-space, parametrised by $t,y,z$. The planar ``slices'', given by fixing $y$ and $z$ in the first two coordinates of (\ref{eq:butterfly-param}),
are curves parametrised by $t$. As $y$ and $z$ vary near the value 0 we can observe how these curves evolve. When $y=z=0$ the curve has a $(4,5)$ singularity
(switching the coordinates), which is the same as the singularity in the central diagram of Figure~\ref{fig:butterfly1}. A simple Desmos demonstration enabling you to
draw the curves for different $(y,z)$,
 appearing in Desmos as $(Y,Z)$, is at \cite{Desmos2}. A different version
\cite{Desmos3} takes  the point $(Y,Z)$ for a tour
around the origin $Y=Z=0$.
See also 
Figure~\ref{fig:desmos-butterfly} for a ``still'' from \cite{Desmos3}.

\begin{figure}[h!]
\begin{center}
\includegraphics[width=1.6in]{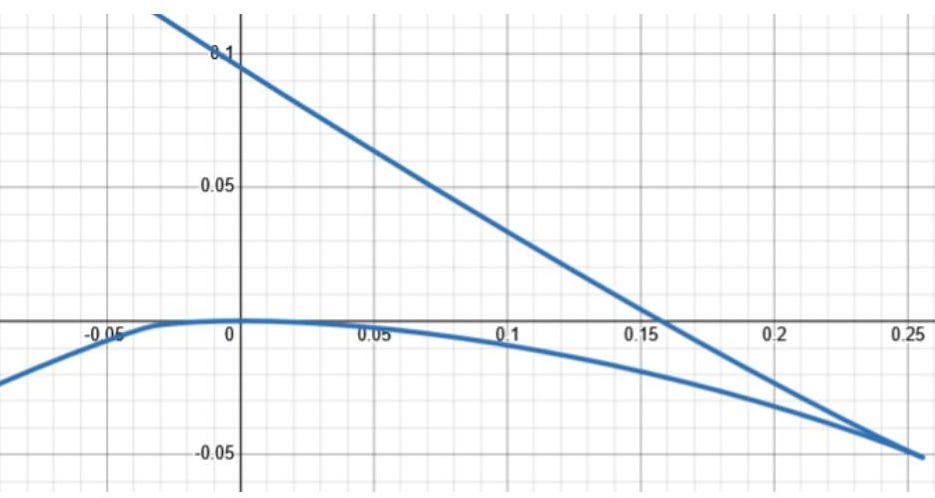}
\includegraphics[width=1.6in]{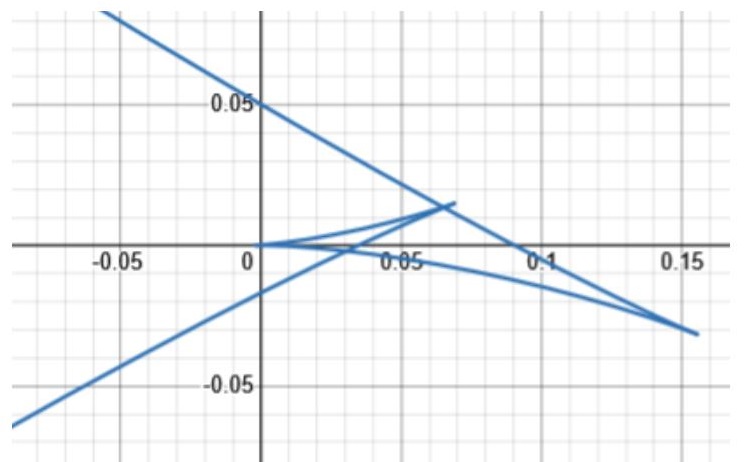}
\includegraphics[width=1.6in]{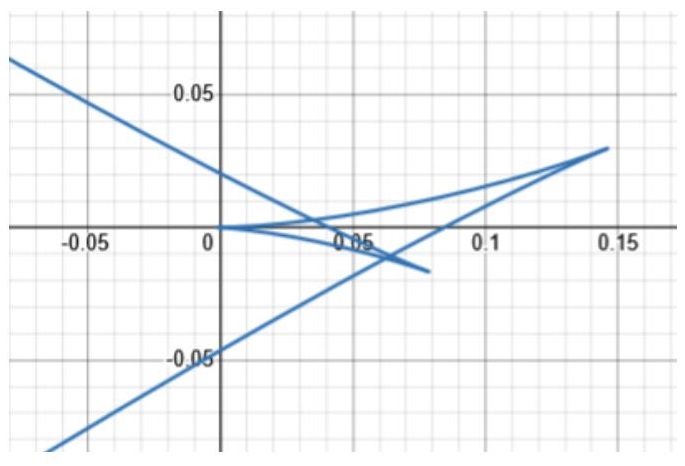}
\includegraphics[width=1.6in]{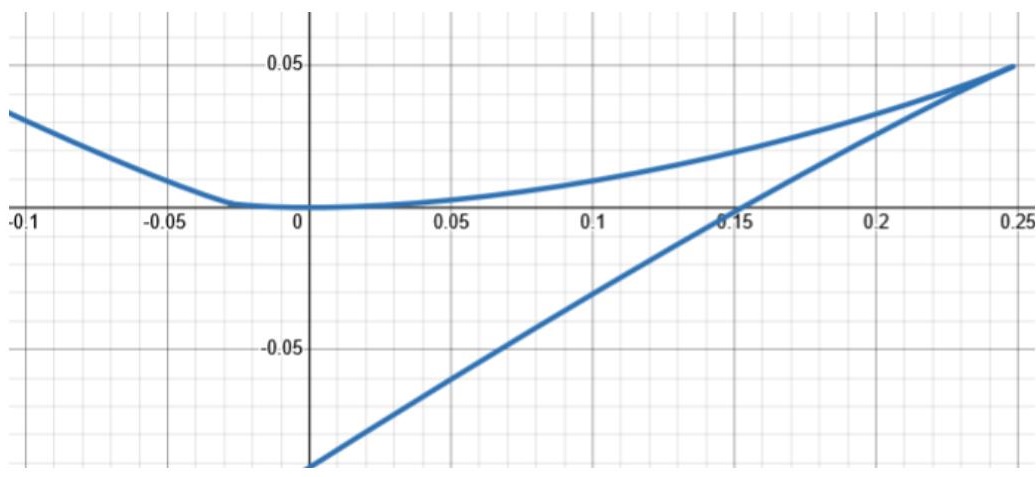}
\end{center}
\vspace{-0.2in}
\caption{\small 
Four diagrams from a tour of the origin $(y, z)=(0.5\cos a ,\ 0.5\sin a)$ for, left to right,
$a=4.3, \ 4.6, \ 4.8, \ 5.1$ exhibiting two swallowtail transitions in a neighborhood of
a butterfly singularity. The symmetric diagram $a=4.712$ 
is in Figure~\ref{fig:desmos-butterfly}.}
\label{fig:butterfly-standard}
\end{figure}

Two swallowtail transitions for the family (\ref{eq:H}), 
observed  in a the circular tour of $(y,z)=(0,0)$, are also displayed in
 Figure~\ref{fig:butterfly-standard}.

We now show how to add a parameter to the 2-circles family so that the full
structure of a more complicated singularity can be unfolded.

\subsection{Adding an extra parameter to the 2-circles construction} \label{ss:unf}

We can extend the function $F(t,x,y,r)$ in (\ref{eq:F-2circles}), which gives the family
of lines joining corresponding points of two  circles both centered at the origin and
of radius 1 and $r$, by allowing
the center of the second circle to move to $(0,d)$. That is, the line now joins
$(\cos t, \sin t)$ to $(r\cos mt, d+r\sin mt)$, with equation $\mathcal{F}=0$ where
\begin{equation}
\mathcal{F}(t,x,y,r,d)=x(r\sin mt -sin t +d)-y(r\cos mt -\cos t)-
r\cos t\sin mt+r\sin t\cos mt -d\cos t.
\end{equation}\label{eq:Ftxyrd} 
There is a correspondingly more complicated formula for the discriminant
$\mathcal F=\mathcal{F}_t=0$. 

The addition of the parameter $d$ can be described as {\em breaking the symmetry}
of the 2-circles envelope construction and it confers an extra degree of freedom on the
envelope curve.

\begin{figure}[H]
\begin{center}
\includegraphics[width=1.1in]{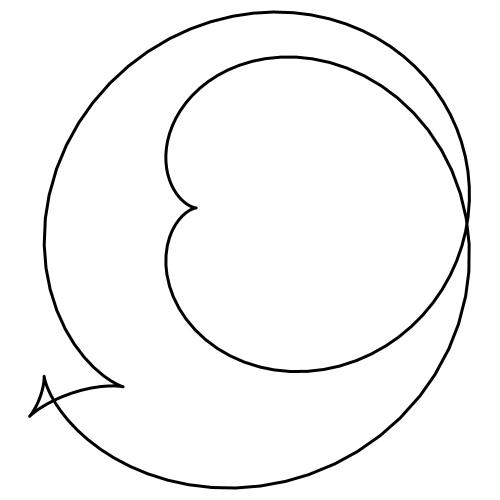}
\includegraphics[width=1.1in]{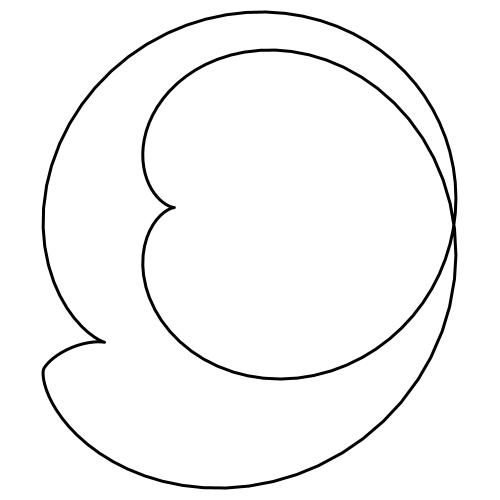}
\includegraphics[width=1.1in]{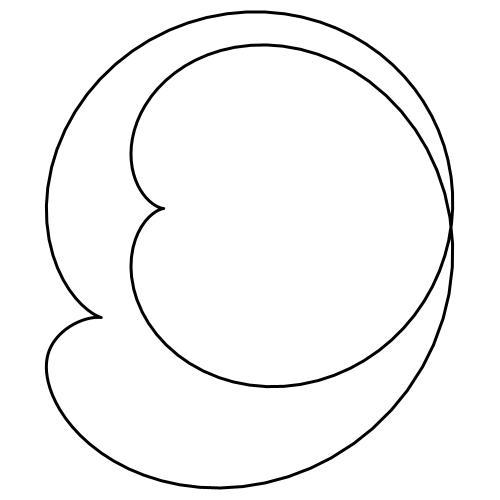}
\end{center}
\hspace{1.7in}$r=2.5,d=0.5$ \ \   $r=3,d=0.5$ \ \   $r=3.5, d=0.5$
\begin{center}
\includegraphics[width=1.1in]{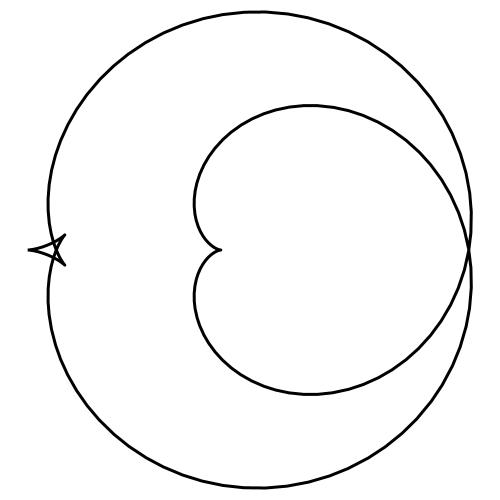}
\includegraphics[width=1.1in]{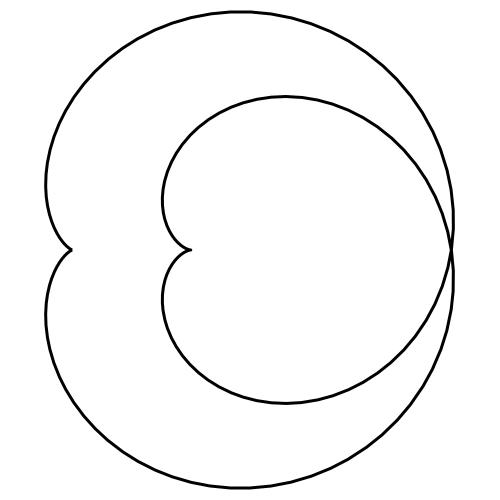}
\includegraphics[width=1.1in]{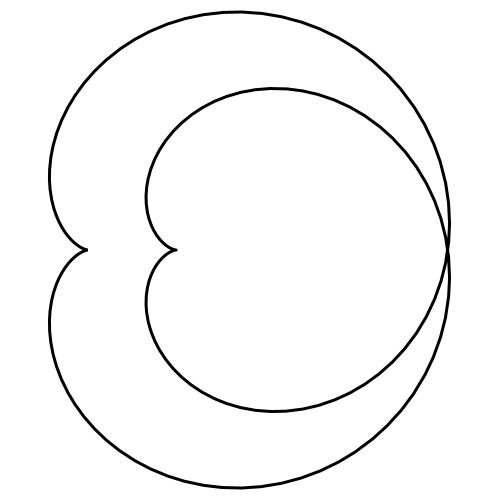}
\end{center}
\hspace{1.7in}$r=2,d=0$ \ \ \ \  \ \  $r=2.5,d=0$ \ \  \ \  \ \ $r=3, d=0$
\begin{center}
\includegraphics[width=1.1in]{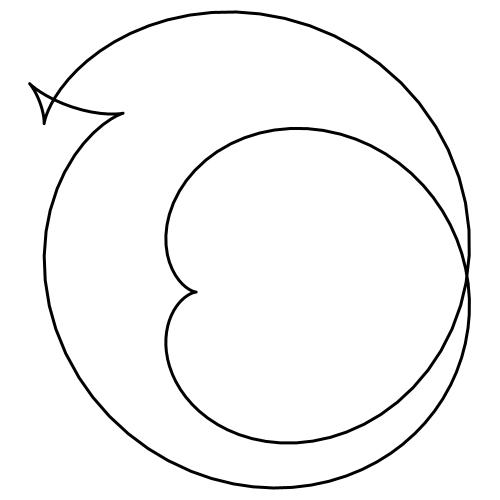}
\includegraphics[width=1.1in]{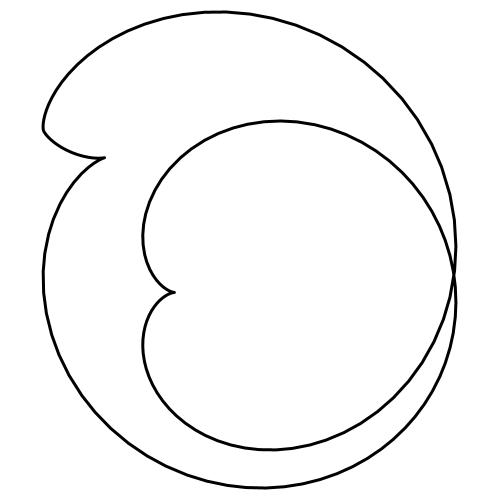}
\includegraphics[width=1.1in]{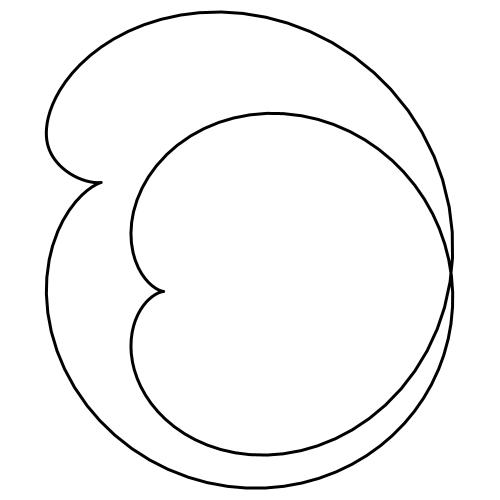}
\end{center}
\hspace{1.6in}$r=2.5,d=-0.5$ \    $r=3,d=-0.5$ \    $r=3, d=-0.5$
\caption{\small A tour of the origin in the $(r,d)$-plane, showing the envelope obtained from  the unit circle and the circle radius $r$ near to 2.5 and center $(0,d)$
near to $(0,0)$.  Throughout,
$m=3/4$ and all figures are scaled 
to be approximately the same size. The diagram as a whole is called a ``clock diagram''
for the butterfly singularity which is to the left in the central diagram. It shows a ``universal
unfolding'' of that singularity, revealing all ``nearby singularities'' to the butterfly, that is
singularities which can occur in any small perturbation.  Swallowtail transitions occur top left to top right, and bottom left to
bottom right.}
\label{fig:butterfly1}
\end{figure}

 Figure~\ref{fig:butterfly1} displays a ``clock diagram'' of plane curves obtained
  by taking $m=3/4$ and taking values of $r$ and $d$ 
  close to $(r,d)=(2.5,0)$. This base value makes the singularity of the
  envelope at $(-5.7,0)$ of type (4,5): see (\ref{eq:taylor}).
  A two-dimensional array of figures is needed to display
the independent variation of the two parameters $r,d$.

\subsection{Some more technical comments}\label{ss:tech}

How can we disinguish between a family such as $H$ in (\ref{eq:H}), 
 
Here we shall indicate how it can be decided whether an unfolding will display
 ``all nearby singularities''  to a given starting singularity, that is, reveal the full
 inner structure of that singularity. The methods are simple but the justification is more
 technical and we must refer to books such as 
 \cite[p.148--151]{CS},\cite{Mon} for the details.
 There are also some details on the webpage \cite{unf} though these refer to
 a slightly different classification.

Using the methods of \cite[Ch.7]{CS} it can be shown
that the two parameters $r$ and $d$ are always sufficient to produce the
full unfolding of the butterfly singularity exhibited in this figure. In this diagram there 
are two instances of swallowtail transitions: top left to top right and bottom
left to bottom right. Two cusps disappear via a swallowtail point to a smooth piece
of curve.

Besides books such a \cite{BL,Mon,PS} there is  much
information and good illustrations of these singularities, or
``catastrophes'' on the internet, such as~\cite{W}. 

\medskip\noindent
{\sc Acknowledgements}  \ This article is based on a Work Experience placement 
undertaken by Wettig, supervised by Giblin, at Liverpool
University. Subsequently Giblin has added more material, especially in the
last sections of the article. Wettig is  grateful to his High School, Liverpool College, for suggesting that he undertake a work placement at Liverpool University, and  both authors
acknowledge the University for the use of computing and other facilities.

\medskip\noindent
Peter Giblin, Department of Mathematical Sciences, The University of Liverpool, U.K. pjgiblin@liv.ac.uk \\
Alexander Wettig, Department of Computer Science, Princeton University, U.S.A.  awettig@cs.princeton.edu
\end{document}